\documentstyle[11pt]{article}
\newcommand{\doublespace}{
   \renewcommand{\baselinestretch}{1.2}
   \large\normalsize}

\def \Z{\Bbb Z}
\def \C{\Bbb C}

\def \Q{\Bbb Q}

\def \M{\cal M}

\def \wt{{\rm wt}}

\def \Res{{\rm Res}}

\def \End{{\rm End}}
\def \Hom{{\rm Hom}}

\def \Ind {{\rm Ind}}

\def \Aut{{\rm Aut}}

\def \<{\langle}
\def \>{\rangle}

\def \a{\alpha }

\def \l{\lambda }
\def \L{\Lambda }
\def \g{\gamma}
\def \b{\beta }

\def \p{\phi}

\def \qed{\mbox{ $\square$}}
\def \pf {\noindent {\bf Proof:} \,}
\def \cg{\chi_g}
\def \cg'{\chi'_g}

\def \o{\otimes}
\def \d{\delta}

\input amssym.def
\input amssym
\doublespace
\begin{document}
\newtheorem{th}{Theorem}[section]
\newtheorem{thmn}{Theorem}
\newtheorem{prop}[th]{Proposition}
\newtheorem{cor}[th]{Corollary}
\newtheorem{lem}[th]{Lemma}
\newtheorem{rem}[th]{Remark}
\newtheorem{de}[th]{Definition}
\newtheorem{hy}[th]{Hypothesis}
\begin{center}
{\Large {\bf The relationship between skew group algebras and orbifold theory}
} \\
\vspace{0.5cm}
 Gaywalee Yamskulna\footnote{Supported by The Development and Promotion of Science and Talents Project Fellowship from the government of Thailand.}
\\
Mathematical Sciences Research Institute, Berkeley, CA 94720
\end{center}

\hspace{1.5 cm}
\begin{abstract}

Let $V$ be a simple vertex operator algebra and let $G$ be a finite automorphism group of $V$.
In [DY], it was shown that any irreducible $V$-module is a completely reducible $V^G$-module where $V^G$ is the $G$-invariant sub-vertex operator algebra of $V$.
In this paper, we give an alternative proof of this fact using the theory of skew group algebras.
We also extend this result to any irreducible $g$-twisted $V$-module when $g$ is in the center of $G$ and $V$ is a $g$-rational vertex operator algebra.
\end{abstract}
\section{Introduction}

In orbifold theory, one studies a simple vertex operator algebra (VOA) $V$, a finite automorphism group $G$ and the $G$-invariant sub-vertex operator algebra $V^G$.
One of the important problems in this area is to determine the module category for $V^G$.
In considering this question, the twisted $V$-modules become important because they restrict to $V^G$-modules.
Roughly speaking, for $g\in G$, a $g$-twisted $V$-module $M$ is a $\C$-graded infinite-dimensional vector space such that each element $v\in V$ acts on $M$ by $Y_M(v,z)=\sum_{n\in \Q}v_nz^{-n-1}$, where $z$ is a formal variable
and $v_n\in \End M$. If $g$ is the identity element, then $g$-twisted $V$-modules are called $V$-modules.

A vertex operator algebra $V$ is called rational if any $V$-module is completely reducible.
The {\em main conjecture} in orbifold theory is the following: {\em If $V$ is rational then $V^G$ is rational}.
{\em Moreover, every irreducible $V^G$-module is contained in some irreducible $g$-twisted $V$-module for some $g\in G$ (cf. [DVVV])}.
This suggests that for each $g\in G$ one must study each irreducible $g$-twisted $V$-module as a $V^G$-module. In order to approach this conjecture, it is natural to start with the case when the irreducible $g$-twisted $V$-module is the simple vertex operator algebra $V$, itself.
This case was done by Dong, Li and Mason in [DLM1].
For a given simple VOA $V$ and a finite-dimensional compact-Lie group $G$ such that the action of $G$ on $V$ is continuous,
they showed that $V$ is a completely reducible $V^G$-module. To be more precise, they
showed that every simple $G$-module occurs in $V$ and the multiplicity
space of each simple $G$-module in $V$ is an irreducible
$V^G$-module. Furthermore, they showed that inequivalent irreducible $G$-modules
produce inequivalent $V^G$-modules.
In [DY], Dong and the author extended the results in
[DLM1] to any irreducible $V$-module $M$ when $G$ is a finite
automorphism group.
Several duality theorems of Schur-Weyl type were also obtained.
In particular, all inequivalent irreducible $V^G$-modules which occur as $V^G$-submodules of $M$ were classified.
Moreover, we established relationships between the irreducible $V^G$-submodules of $M$ and the irreducible $V^G$-submodules of another irreducible $V$-module $N$. Indeed, the results in [DY] hold for any irreducible $g$-twisted $V$-modules when $g$ is in the center of $G$ (cf. [Y]).

In this paper, we use the theory of skew group algebras to provide
an alternative way of showing that any irreducible $g$-twisted
$V$-module is a completely reducible $V^G$-module when $g$ is in
the center of $G$ and $V$ is a $g$-rational simple vertex operator
algebra. Let $A$ be a finite-dimensional algebra over $\C$ and $G$
be a finite automorphism group of $A$. To these objects we
associate the skew group algebra $A\rtimes G$ as defined in
section 2. We then consider a finite-dimensional simple $A$-module
$N$ and its inertia subgroup $G_N=\{h\in G\ \ |\ \ ^hN\cong N\ \
\mbox{as}\ \ A\mbox{-modules}\}$. Here, $^hN$ is a new simple
$A$-module associated to $N$ and $h$. The group $G_N$ acts
projectively on $N$. We let $\a$ be the corresponding 2-cocycle
whose values are roots of unity and we denote by $I$ the set of
inequivalent simple $\C^{\a}[G_N]$-modules occurring in $N$. For
each $W\in I$ ,  we let $N_W$ be the multiplicity space of $W$ in
$N$. The space $N_W$ is, in fact, an $A^G$-module (see Lemma 3.1)
where, $A^G$ is the $G$-invariant sub-algebra of $A$. By using the
invariant theory of skew group algebras, we prove the following
result.
\begin{th}\label{th11} If $A$ is semisimple, then the space $N_W$ is a simple $A^G$-module for all $W\in I$.
\end{th}

In [DLM4], Dong, Li, and Mason constructed a series of associative algebras $A_{g,n}(V)$ from a given vertex operator algebra $V$. These associative algebras are quotient spaces of $V$ associated with an element $g\in G$ and $n\in {1\over T}\Z_+$
where $T$ is the order of $g$ and $\Z_+$ is the set of nonnegative integers.
These algebras were first introduce by Zhu for the case when $g$ is the identity element and $n=0$ (cf. [Z]).
The important properties of the $A_{g,n}(V)$ that are used in this paper are the following: let $M=\bigoplus_{n\in {1\over T}\Z_+} M(n)$ be an admissible $g$-twisted $V$-module with $M(0)\neq 0$. Then

1) $M$ is an irreducible $g$-twisted $V$-module if and only if each $M(n)$ is a simple $A_{g,n}(V)$-module for all $n\in {1\over T}\Z_+$.

2) $M^1$ and $M^2$ are isomorphic as irreducible $g$-twisted
$V$-modules if and only if $M^1(n)$ and $M^2(n)$ are isomorphic as
nonzero simple $A_{g,n}(V)$-modules.

3) If $V$ is a $g$-rational VOA, then $A_{g,n}(V)$ is a finite-dimensional semisimple associative algebra.

This key results allow us to reduce an infinite dimensional problem to a finite dimensional one.

Now we assume that $g$ is the element in the center of $G$ and
that $V$ is a $g$-rational simple vertex operator algebra. Set
$G_M=\{h\in G\ \ |\ \ M\circ h\cong M\ \ \mbox{as}\ \
g\mbox{-twisted}\ \ V\mbox{-modules}\}$. Here, $M\circ h$ is a new
irreducible $g$-twisted $V$-module associated with $g$ and $M$.
The group $G_M$ acts projectively on $M$. We let $\a_M$ be the
corresponding 2-cocycle whose values are roots of unity and denote
by $\L_M$ the set of all inequivalent simple
$\C^{\a_M}[G_M]$-modules. Hence $M=\bigoplus_{W\in \L_M}M^W$ where
$M^W$ is the sum of all simple $\C^{\a_M}[G_M]$-submodules of $M$
isomorphic to $W$. Furthermore, $M\cong \bigoplus_{W\in \L_M}W\o
M_W$ as $\C^{\a_M}[G_M]\o V^G$-modules where $M_W$ be the
multiplicity space of $W$ in $M$. By using the properties of the
associative algebras $A_{g,n}(V)$ combined with Theorem
\ref{th11}, we obtain the following.
\begin{th}\label{th12} The space $M_W$ is an irreducible $V^G$-module.
\end{th}
This immediately implies that
\begin{th}\label{th13}
If $g$ is in the center of $G$ and $V$ is a $g$-rational
simple VOA, then every irreducible $g$-twisted $V$-module is a
completely reducible $V^G$-module.
\end{th}

This paper is organized as follows. In the second section, we recall the definition of a skew group
algebra $A\rtimes G$. We also discuss the relation between the
modules of $A\rtimes G$ and $A$-modules, and we review several
important theorems of $A^G$-modules. In the third section, we
prove Theorem \ref{th11}. In the fourth section, we recall the
definitions of vertex operator algebras and their twisted modules,
and we discuss the representation theory of VOAs. We also review
the construction of the series of associative algebra $A_{g,n}(V)$
and their important properties. In the last section, Theorem
\ref{th11} is used to prove Theorem \ref{th12} and Theorem
\ref{th13}.

{\bf Acknowledgments:} I would like to express my gratitude to my advisor Professor Chongying Dong for his helpful advice and comments. I would also like to thank Professor Georgia Benkart for providing me with a very useful preprint of her text book. As well, I am grateful to Professor Susan Montgomery for her advice about fixed ring theory.

\section{Skew Group Algebras}

Let $A$ be an algebra and let $G$ be a finite automorphism group of $A$.
We give the definition of the {\em skew group algebra} $A\rtimes G$.
We also discuss the relations between modules of $A\rtimes G$ and $A$-modules.
Then, we study the $G$-invariant sub-algebra of $A$ and its modules.
\begin{de} {[BR, Mo]}
Let $A$ be an algebra over a field $F$. Let $G$ be a finite
automorphism group of $A$. The {\em skew group algebra} $A\rtimes
G$ is the set $\{\ \ \sum_{g\in G}a_g g\ \ |\ \ a_g\in A\ \ \}$ of
all finite $A$-linear combinations of elements of $G$ with
multiplication that is the linear extension of
$$agbh=ag(b)gh$$
for $a, b\in A$ and $g,h\in G.$
\end{de}

For an $A$-module $M$ and $g\in G$, we set $^gM=M$ as vector
spaces over $F$. We define the $A$-action on $^gM$ in the
following way : for $a\in A$, $m\in\ \ ^gM$,
\begin{equation}\label{sk2}
 a*_gm=g^{-1}(a)m.
\end{equation}

Now assume that $M$ is a finite-dimensional {\em simple} $A$-module and the
field $F$ is algebraically closed.
 The {\em{inertia subgroup}} of $M$ is $ G_M=\{ \ \ h\in G\ \ |\ \ ^hM\cong M\ \
\mbox{as} \ \ A\mbox{-modules} \ \ \}.$ By Schur's Lemma, we
conclude that $\Hom_A(^hM,M)$ is one-dimensional. For each $h\in
G_M$, we fix an isomorphism $\phi (h):\ \ ^hM\rightarrow M$. Here,
$\p(h)$ satisfies $\p(h)(h^{-1}(a)m)=\p(h)(a*_hm)=a\p(h)m.$ This
implies that
\begin{equation}\label{sk1}
\phi(h)a=h(a)\phi(h).
\end{equation}
Any linear transformation $\psi$ on $M$ satisfying $\psi a
=h(a)\psi$ for all $a\in A$ must belong to $\Hom_A(^hM, M)$.
Hence, $\psi$ is a multiple of $\p (h)$. Since $\p(h)\p(k)a =
(hk)(a)\p(h)\p(k),$ we see that $\p(h)\p(k)\in \Hom_A(^{hk}M,M)$.
Furthermore, there is $\a(h,k)\in F^{\times}$ such that
$$ \p(h)\p(k)=\a(h,k)\p(hk).$$
By associativity of composition, we have
$$\a(h,k)\a(hk,l)=\a(h,kl)\a(k,l)\ \ \mbox{for\ \ all} \ \ h,k,l\in
G_M. $$ Note that one can choose $\p(1)=1_M$ in order to obtain
$\a(1,h)=\a(h,1)=1$ for $h\in G$. Thus $F[G_M]$ acts projectively
on $M$ and $\a$ is a 2-cocycle. By using this $\a$ we define a
twisted group algebra $F^{\a^{-1}}[G_M]$ over $F$ with basis $c_h,
h\in G_M$, and with multiplication defined by
$$ c_h c_k= \a(h,k)^{-1}c_{hk}\ \ \ \ \ \ \ h,k \in G_M.$$

\begin{prop}{[BR]}\label{pr3.3}
Let $M$ be a simple $A$-module.
For any $F^{\a^{-1}}[G_M]$-module $V$, there is an $(A\rtimes G_M)$-module action
on $M\otimes V$ defined by
\begin{equation}\label{br1}
(ah)(m\otimes v)=a\p(h)(m)\otimes c_hv.
\end{equation}
\end{prop}
\pf Let $a,b\in A$ and $h,k\in G_M$. We have
\begin{eqnarray*}
(ah)((bk)(m\o v))&=&(ah)(b\p(k)(m)\o c_kv)\\
&=&a\p(h)(b\p(k)(m))\o c_h c_k v\\
&=&ah(b)\p(h)\p(k)(m)\o \a(h,k)^{-1}c_{hk}v\\
&=&ah(b)\a(h,k)\p(hk)\o \a(h,k)^{-1}c_{hk}v\\
&=&ah(b)\p(hk)(m)\o c_{hk}v\\
&=&ah(b)hk(m\o v)\\
&=&(ahbk)(m\o v).\ \ \qed
\end{eqnarray*}

Next we discuss the modules of the algebra $A\rtimes G$ and their relation to $A$-modules.
In fact, we recall that simple $A\rtimes G$-modules can be obtained from the simple $A$-modules and certain twisted group algebras.
\begin{th}{[BR]}
Assume $N$ is a finite-dimensional simple $A\rtimes G$-module. Let
$A^{\l}$ be a simple $A$-submodule of $N$, and let $H$ denotes the
inertia subgroup of $A^{\l}$. Suppose that the isomorphisms
$^hA^{\l}\stackrel{\p(h)}{\rightarrow} A^{\l}, h\in H$, determine
the cocycle $\a$ by $\p(h)\p(k)=\a(h,k)\p(hk)$. Then there is a
simple $F^{\a^{-1}}[H]$-module $H^{\nu}$ so that
$$N \cong \Ind_{A\rtimes H}^{A\rtimes G}(A^{\l} \otimes H^{\nu}).$$
\end{th}

\begin{cor}
In fact, $H^{\nu}=\Hom_A(A^{\l},P)$ where $P=\sum_{h\in H}hA^{\l}$.
\end{cor}

Next, we discuss the $G$-invariants, $A^G=\{a\in A\ \ |\ \ g(a)=a
\mbox{ for\ \ all\ \ } g\in G\}$. We \em{assume that $|G|^{-1}\in F$.}
Set
$$ e=|G|^{-1}\sum_{g\in G}g \in A\rtimes G.$$
Clearly, $e$ is an idempotent of $A\rtimes G$. Furthermore,
$ea=ae$ and $eh=he$ for all $a\in A^G$, $h\in G$. If $a$ is an
arbitrary element of $A$, then $\sum_{g\in G}g(a)\in A^G$.
Moreover, we obtain the following.
\begin{prop}{[BR]}
\ \

i) The map $\Phi : A^G\rightarrow e(A\rtimes G)e,\ \ a\mapsto ae$
is an algebra isomorphism.

ii) The map $\Psi: A\rightarrow (A\rtimes G)e,\ \ a\mapsto ae$, is
an isomorphism of $(A\rtimes G, A^G)$-bimodules.
\end{prop}

The following theorem is the key theorem of this paper. This
theorem can be found in [BR].
\begin{th} \label{thsi}
( {\em Invariant Theory} ) Let $M$ be a simple module for $A\rtimes
G$.
Then $eM$ is a simple module for $e(A\rtimes G)e$.
Furthermore, every simple $e(A\rtimes G)e$-module is obtained in this
way.
\end{th}

\section{The Main Theorem}
For the rest of this paper, we set $F=\C$. We also assume that $A$
{\em is semisimple}. Let $M$ be a {\em finite-dimensional simple}
$A$-module and let $G_M$ be its inertia subgroup. By the previous
section, we know that $G_M$ acts projectively on $M$ and we let
$\a_M$ be the corresponding 2-cocycle whose values are roots of unity.
Since $\C^{\a_M}[G_M]$ is semisimple, $M$ is a semisimple $\C^{\a_M}[G_M]$-module. Let
$J_{M,\a _M}$ be the set of irreducible $\a_M$-characters
$\lambda$ of $\C^{\a_M}[G_M]$ such that the corresponding
simple-modules $W_{\l}$ occur in $M$. We write $M=\oplus_{\l\in
J_{M,\a _M} }M^{\l}$ where $M^{\l}$ is the sum of simple
$\C^{\a_M}[G_M]$-submodules of $M$ isomorphic to $W_{\l}$. Let
$M_{\l}$ be the multiplicity space of $W_{\l}$ in $M$ (ie.
$M_{\l}=\Hom_{\C^{\a_M}[G_M]}(W_{\l}, M)$). We can realize
$M_{\l}$ as a subspace of $M$ in the following way. Let $w\in
W_{\l}$ be a fixed nonzero vector. We identify
$\Hom_{\C^{\a_M}[G_M]}(W_{\l}, M)$ with the subspace $\{f(w) |f\in
\Hom_{\C^{\a_M}[G_M]}(W_{\l},M)\}$ of $M^{\l}$.
\begin{lem}
$M_{\l}$ is an $A^G$-module.
\end{lem}
\pf Since $\p(h)a=a\p(h)$ for all $a\in A^G$ (cf. (\ref{sk1})), hence the actions of $\C^{\a_M}[G_M]$ and $A^G$ on $M$ are commutative and the lemma follows immediately. \qed
\begin{rem}
$M_{\l}$ is also an $A^{G_M}$-module.
\end{rem}

Observe that $M^{\l}$ and $M_{\l}\o W_{\l}$ are isomorphic as both
$\C^{\a_M}[G_M]$-modules and $A^G$-modules. Here, $A^G$ acts on the
first tensor factor of $M_{\l}\o W_{\l}$ and $\C^{\a_M}[G_M]$
acts on the second tensor factor. We now {\em identify $M^{\l}$
with $M_{\l}\o W_{\l}$} as $A^G\o \C^{\a_M}[G_M]$-modules.

Let $\{\bar{g}|g\in G\}$ be the basis for $\C^{\a}[G]$-module.
Recall that for a $\C^{\a}[G]$-module $V$, $V^*=\Hom_{\C}(V, \C)$
becomes a $\C^{\a^{-1}}[G]$-module via
$$(c_g \psi)(v)=\psi(\bar{g}^{-1}v)$$ for all $g\in G, v\in V,
\psi\in V^*$. Here, $\{c_g|g\in G\}$ is a basis of
$\C^{\a^{-1}}[G]$. The dual space $V^*$ is called the
\em{contragredient module} of $V$.

For a fixed $\g\in J_{M,\a_M}$, we let $\g^*$ be the
${\a_M}^{-1}$-character of $\C^{{\a_M}^{-1}}[G_M]$ dual to $\g$.
The corresponding $\C^{{\a_M}^{-1}}[G_M]$-module is
${W_{\g}}^*$. By Proposition (\ref{pr3.3}),
we conclude that $M\o W_{\g}^*$ is an $A\rtimes G_M$-module.

\begin{th} If $A$ is a semisimple associative algebra, then $\Ind_{A\rtimes G_M}^{A\rtimes G}M\o
W_{\g}^*$ is a simple $A\rtimes G$-module.
\end{th}

\pf For any nonzero $w\in W_{\g}^*$, we let $M_h=M\o c_h w$ for
$h\in G_M$. Then we have $M\o W_{\g}^* =\sum_{h\in G_M}M_h.$ Note
that for $h\in G_M$, $M_h$ is isomorphic to $M$ as $A$-modules.
Let $G=\bigcup_{i=1}^kg_iG_M$ be a left coset decomposition.
Hence, $$\Ind_{A\rtimes G_M}^{A\rtimes G}M\o
W_{\g}^*=\oplus_{i=1}^k g_i\o \sum_{h\in
G_M}M_h=\oplus_{i=1}^k\sum_{h\in G_M}g_i\o M_h.$$ The space $g_i\o
M_h$ is an $A$-module under the following action: for $a\in A$,
$g_i\o m\o c_h w \in g_i \o M_h$,
$$a\cdot(g_i\o m\o c_h w)=g_i\o
g_i^{-1}(a)m\o c_hw.$$ Indeed, $g_i\o M_h$ is isomorphic to $^{g_i}M$ as $A$-modules.
Hence, $g_i\o M_h$ is a simple $A$-module.

Finally, we will show that $\Ind_{A\rtimes G_M}^{A\rtimes G}M\o
W_{\g}^*$ is a simple $A\rtimes G$-module. We let $N$ be an
$A\rtimes G$-submodule of $Ind_{A\rtimes G_M}^{A\rtimes G}M\o
W_{\g}^*$.
Consequently, $N$ is an $A$-module. Since $A$ is semisimple,
$N$ contains a simple $A$-submodule of $\oplus_{i=1}^k\sum_{h\in
G_M}g_i\o M_h$. Hence, there exist $j\in \{1,...,k\}$ and $ h'\in G_M$ such that $g_j \o M_{h'}$
is contained in $N$. Then  $1\o M_1=(g_jh')^{-1}\cdot g_j \o M_{h'}\subset N$.
These imply that $\oplus_{i=1}^k\sum_{h\in G_M}g_i\o M_h \subset N$. Therefore, $Ind_{A\rtimes G_M}^{A\rtimes G}M\o
W_{\g}^*$ is a simple $A\rtimes G$-module. \qed

For any finite group $G$, and any $G$-module $W$, we denote the
$G$-invariant submodule of $W$ by $Inv_G(W)$ (i.e, $Inv_G(W)=\{ x\in
W\ \ | \ \ gx=x \mbox{\ \ for\ \ all\ \ }g\in G\}$).

\begin{prop}\label{l44}{[K]} For any finite dimensional $\C^{\a}[G]$-modules $N$, $M$, we have the natural isomorphisms
$$\Hom_{\C^{\a}[G]}(M,N)\cong \Hom_{\C[G]}(\C,M^*\o N)\cong Inv_G(M^*\o N).$$
Moreover, $Inv_G(M^*\o N)=\left(\sum_{g\in G}g\right)\cdot (M^*\o N).$
\end{prop}

Set $${\mathcal M }=\Ind_{A\rtimes G_M}^{A\rtimes G}M\o
W_{\g}^*.$$ Then
$${\mathcal M}=\bigoplus_{\l\in J_{M,\a _M}}\bigoplus_{i=1}^k g_i\o M_{\l}\o W_{\l}\o W_{\g}^*$$ where $G=\bigcup_{i=1}^kg_iG_M$ is a left coset decomposition.
Note that for any $\l\in
J_{M,\a _M}$, $W_{\l}\o W_{\g}^*$ is a $\C[G_M]$-module.

\begin{th} \label{thms}
( The Main Theorem ) For any $\g\in J_{M,\a_M}$, the space $M_{\g}$
is a simple $A^G$-module.
\end{th}
\pf
Let $\sum_{\l\in J_{M,\a _M}}\sum_{i=1}^k g_i\o a_i^{\l}\o
b_i^{\l}\o w_i^{\g}\in\mathcal M. $
Then
\begin{eqnarray*}
e(\sum_{\l\in J_{M,\a _M}}\sum_{i=1}^k g_i\o a_i^{\l}\o b_i^{\l}\o
w_i^{\g})&=&\sum_{\l\in J_{M,\a _M}}\sum_{i=1}^k e\o a_i^{\l}\o
b_i^{\l}\o w_i^{\g}\\
&=&\sum_{\l\in J_{M,\a _M}}\sum_{i=1}^k e^2\o a_i^{\l}\o b_i^{\l}\o
w_i^{\g}\\
&=&\sum_{\l\in J_{M,\a _M}}\sum_{i=1}^k e({1\over{|G|}}\sum_{j=1}^k
g_j)(\sum_{h\in G_M}h) \o a_i^{\l}\o b_i^{\l}\o w_i^{\g}\\
&=&\sum_{\l\in J_{M,\a _M}}\sum_{i=1}^k e({1\over{|G|}}\sum_{j=1}^k
g_j)\o a_i^{\l}\o (\sum_{h\in G_M}h) (b_i^{\l}\o w_i^{\g})\\
&=&\sum_{\l\in J_{M,\a _M}}\sum_{i=1}^k {1\over{|G_M|}}e\o a_i^{\l}\o
(\sum_{h\in G_M}h) (b_i^{\l}\o w_i^{\g})\\
&=&\sum_{i=1}^k e\o a_i^{\g}\o{1\over|G_M|}(\sum_{h\in G_M}h) (b_i^{\g}\o
w_i^{\g})\ \ (\mbox{by\ \ Proposition \ref{l44}}).\\
\end{eqnarray*}
We have $e({\mathcal M})\subset e\o M_{\g}\o Inv_G(W_{\g}\o
W_{\g}^*)$. In fact, $e({\mathcal M})=e\o M_{\g}\o Inv_G(W_{\g}\o
W_{\g}^*)$. Since $\M$ is a simple $A\rtimes G$-module, Theorem
\ref{thsi} implies that $e\o M_{\g}\o Inv_G(W_{\g}\o W_{\g}^*)$ is
a simple $A^G$-module. Hence $M_{\g}$ is a simple $A^G$-module.
\qed

\begin{cor}
$M$ is a completely reducible $A^G$-module.
\end{cor}
\section{Vertex operator algebras and related results}

We begin by reviewing the definitions of vertex operator algebras and their twisted modules.
We also recall some relevant results about the representation theory of vertex operator algebras.
Furthermore, for a fixed VOA $V$ we review the construction of a series of associative algebra $A_{g,n}(V)$ and some of their main properties which play an important role in the next section.

For a vector space $W$, we let $W[[z, z^{-1}]]$ be the space of
$W$-valued formal series in arbitrary integral powers of $z$.
\begin{de} {[FLM, FHL]} A {\it vertex operator algebra} is a
${\Z}$-graded vector
space
$$V=\oplus_{n\in{\Z}}V_n;$$ such that
\begin{eqnarray}
\dim\,V_n &<&\infty\ \ \mbox{ and}\\
V_n&=& 0\ \ \mbox{ if}\ \ n\ \ \mbox{is\ \ sufficiently\ \ small.}
\end{eqnarray}
Moreover, there is a linear map
\begin{equation}
\begin{array}{l}
V \to (\mbox{End}\,V)[[z,z^{-1}]]\\
v\mapsto Y(v,z)=\displaystyle{\sum_{n\in{\Z}}v_nz^{-n-1}}\ \ \ \
(v_n\in\mbox{End}\,V)
\end{array}
\end{equation}
and two distinguished vectors ${\bf 1}\in V_0$, $\omega \in V_2$
satisfying the following conditions for all $u, v \in V$:
\begin{eqnarray}
u_nv &=&0\ \ \ \ \ \mbox{for}\ \  n\ \ \mbox{sufficiently
large};\label{e2.2}\\
Y({\bf 1},z)&=&1_V;\label{e2.3}\\
Y(v,z){\bf 1}\in V[[z]]&\mbox{and}&\lim_{z\to 0}Y(v,z){\bf 1}=v;
\end{eqnarray}
\begin{equation}\label{jac}
\begin{array}{c}
\displaystyle{z^{-1}_0\delta\left(\frac{z_1-z_2}{z_0}\right)
Y(u,z_1)Y(v,z_2)-z^{-1}_0\delta\left(\frac{z_2-z_1}{-z_0}\right)
Y(v,z_2)Y(u,z_1)}\\
\displaystyle{=z_2^{-1}\delta
\left(\frac{z_1-z_0}{z_2}\right)
Y(Y(u,z_0)v,z_2)}
\end{array}
\end{equation}
(Jacobi identity) where $\delta(z)=\sum_{n\in {\Z}}z^n$ is
the algebraic formulation of the $\delta$-function at 1, and all
binomial
expressions are to be expanded in nonnegative
integral powers of the second variable;
\begin{equation}\label{e2.6}
[L(m),L(n)]=(m-n)L(m+n)+\frac{1}{12}(m^3-m)\delta_{m+n,0}(\mbox{rank}\,V)
\end{equation}
for $m, n\in {\Z},$ where
\begin{equation}
L(n)=\omega_{n+1}\ \ \ \mbox{for}\ \ \ n\in{\Z}, \ \ \
\mbox{i.e.},\ \ \ Y(\omega,z)=\sum_{n\in{\Z}}L(n)z^{-n-2}
\end{equation}
and
\begin{eqnarray}
& &\mbox{rank}\,V\in {\Q};\\
& &L(0)v=nv=(\mbox{wt}\,v)v \ \ \ \mbox{for}\ \ \ v\in V_n\
(n\in{\Z}); \label{3.40}\\
& &\frac{d}{dz}Y(v,z)=Y(L(-1)v,z). \label{3.41}
\end{eqnarray}
\end{de}
We denote the vertex operator algebra just defined by
$(V,Y,\bf{1},\omega)$ (or briefly, by $V$).
The series $Y(v,z)$ are called {\it vertex operators.}
\begin{rem}
The operators $L(n)$ generate a copy of the Virasoro algebra
represented on $V$ with the central charge rank\,$V.$
\end{rem}

\begin{de}
An {\it automorphism} of  $(V,Y,\bf{1},\omega)$ is a linear
isomorphism $g$: $V\to V$  satisfying
\begin{eqnarray*}
gY(v,z)g^{-1}&=&Y(gv,z),\ \ v\in V,\\
g\bf{1}&=&\bf{1},\\
g{\omega}&=&\omega.
\end{eqnarray*}
\end{de}

Let $\Aut(V)$ denote the group of all automorphisms of $V.$
\begin{rem}
For $g\in \Aut(V)$, $g$ commutes with the component operators $L(n)$ of $\omega,$ and
in particular, $g$ preserves the homogeneous spaces $V_n$ which are
the
eigenspaces for $L(0).$ Consequently, each $V_n$ is a module for
$\Aut(V).$
\end{rem}

For $g$, an  automorphism of the VOA $V$ of order $T$, we denote
the decomposition of $V$ into eigenspaces with respect to the action
of $g$ as $V=\bigoplus_{r=0}^{T-1}V^r$ where $V^r=\{v\in V|gv=e^{2\pi ir/T}v\}$.
For a vector space $W$, we denote the space of $W$-valued formal
series in arbitrary complex powers of $z$ by $W\{z\}$.
\begin{de}
A {\em weak $g$-{\it twisted} $V$-module} $M$ is a vector space
equipped with a linear map
\begin{equation}
\begin{array}{l}
V\to (\mbox{End}\,M)\{z\}\label{map}\\
v\mapsto\displaystyle{ Y_M(v,z)=\sum_{n\in\Q}v_nz^{-n-1}\ \ \
(v_n\in\End M)}
\end{array}
\end{equation}
satisfying  axioms analogous to (\ref{e2.2})-(\ref{e2.3}) and
(\ref{jac}).
To describe these, we let $u\in V^r$, $v\in V$ and $w\in M$. Then
\begin{eqnarray}
& &Y_M(u,z)=\sum_{n\in r/T+\Z}u_nz^{-n-1};
\label{1/2}\\
& &u_nw=0\ \ \ \mbox{for}\ \ \ n\ \ \ \mbox{sufficiently\
large};\label{ds1}\\
& &Y_M({\bf 1},z)=1_M;
\end{eqnarray}
\begin{equation}\label{jacm}
\begin{array}{c}
\displaystyle{z^{-1}_0\delta\left(\frac{z_1-z_2}{z_0}\right)
Y_M(u,z_1)Y_M(v,z_2)-z^{-1}_0\delta\left(\frac{z_2-z_1}{-z_0}\right)
Y_M(v,z_2)Y_M(u,z_1)}\\
\displaystyle{=z_2^{-1}\left(\frac{z_1-z_0}{z_2}\right)^{-r/T}\delta\left(\frac{z_1-z_0}{z_2}\right)

Y_M(Y(u,z_0)v,z_2)}.
\end{array}
\end{equation}
 \end{de}
We denote this module by $(M,Y_M),$ or briefly by $M$.
Equation (\ref{jacm}) is called the {\em twisted Jacobi identity}.
If $g$ is the identity element, this reduces to the definition of a
weak $V$-module and (\ref{jacm}) is the untwisted Jacobi identity .

\begin{de}
Suppose that $(M_{i},Y_{i})$ are two weak $g$-twisted $V$-modules,
$i$=1,2.
A homomorphism from $M_{1}$ to $M_{2}$ is a linear map
$f$ : $M_{1}\rightarrow M_{2}$ which satisfies
$$fY_{M_{1}}(v,z)=Y_{M_{2}}(v,z)f$$ for all $v\in V$.
We call $f$ an isomorphism if $f$ is also a linear isomorphism.
\end{de}

\begin{de}\label{d1} An (ordinary) {\em $g$-twisted $V$-module} is a
weak $g$-twisted $V$-module $M$ which carries a $\C$-grading induced
by the
spectrum of $L(0)$. Then $$M=\bigoplus_{\lambda\in \C}M_{\lambda}$$ where
$M_{\lambda}=\{w\in M|L(0)w=\lambda w\},$ dim $M_{\lambda}<\infty.$
Moreover, for fixed $\lambda, M_{{n\over T}+{\lambda}}=0 $ for all
small enough integers $n.$
\end{de}

\begin{lem}\label{l31}
If $M$ is a simple (ordinary) $g$-twisted $V$-module, then $M=\bigoplus_{n=0}^{\infty}M_{\l+{n\over T}}$
for some $\l\in \C$ such that $M_{\l}\neq 0$ and $M_{\l+\frac{n}{T}}=0$ for $n<0$.
\end{lem}

\pf Let $v\in V^r, m\in \frac{r}{T}+\Z$. We recall that, $[L(0),
v_m]=(\wt v -m-1)v_m$ on $M$. By using this relation, we can show
that $M(\b)=\sum_{n\in \Z}M_{\b+\frac{n}{T}}$ is a $g$-twisted
$V$-submodule of $M$. Here, $\b\in \C$. Hence, $M=M(\b)$. By using
the fact that $M_{\b+\frac{n}{T}}=0$ for all small enough integer
$n$. We can choose $\l\in \C$ such that
$M=\bigoplus_{n=0}^{\infty}M_{\l+{n\over T}}$, $M_{\l}\neq 0$, and
$M_{\l+\frac{n}{T}}=0$ for $n<0$. \qed

\begin{de}\label{d2} An {\em admissible $g$-twisted $V$-module} is a
weak $g$-twisted $V$-module $M$ which carries a ${1\over T}\Z_{+}$
grading $M=\bigoplus_{n\in {1\over T}\Z_{+}}M(n)$ satisfying the following
condition:$$v_{m}M(n)\subset M(n+\wt v-m-1)$$ for homogeneous $v\in
V.$ Here, $\Z_+$ is the set of nonnegative integers.
\end{de}
The notion of admissible $g$-twisted $V$-module here is equivalent to the notion of
a module
in [Z] when $g$ is the identity element. Using a grading shift we can always arrange the
grading on $M$ so that $M(0)\ne 0.$ This shift is important
in the study of the algebra $A_{g,n}(V)$ below.

\begin{lem}{[DLM2]}
Any $g$-twisted $V$-module is an admissible $g$-twisted $V$-module.
\end{lem}
So, there is a natural identification of the category
of $g$-twisted $V$-modules with a sub-category of the category
of admissible $g$-twisted $V$-modules.

\begin{de}\label{ration} $V$ is called $g$-rational if every
admissible $g$-twisted $V$-module
is a direct sum of irreducible admissible $g$-twisted $V$-modules.
\end{de}

\begin{th}\label{th3r}{[DLM2]}
If $V$ is $g$-rational then there are only
finitely many inequivalent irreducible admissible $g$-twisted
$V$-modules.
Moreover, each irreducible admissible $g$-twisted $V$-module is an
ordinary $g$-twisted $V$-module.
\end{th}

\begin{prop}\label{4.19}{[L, DM1]} If $M$ is a simple weak $g$-twisted
$V-$module then
$M$ is spanned by $\{u_nm|u\in V,n\in \Q\}$ where $m\in M$ is a fixed
nonzero
vector.
\end {prop}

Next, we recall {\em the construction of the associative
algebra $A_{g,n}(V)$ and some related results}.
This algebra was first introduced by Zhu for the case
$g$ is the identity element and $n$ is zero.
It was later done for the general case by Dong, Li, and Mason.

\setcounter{equation}{0}
Let $V=(V,Y,{\bf 1},\omega)$ be a vertex operator algebra
and $g$ be an automorphism of $V$ of order $T.$ Then $V$ is a direct
sum of eigenspaces of $g:$ $V=\bigoplus_{r=0}^{T-1}V^r$ where $V^r=\{v\in V|gv=e^{2\pi ir/T}v\}.$
Fix $n=l+\frac{i}{T}\in\frac{1}{T}\Z$ with $l$ a nonnegative
integer and $0\leq i\leq T-1.$ For $0\leq r\leq T-1$ we define
\[\d _i(r)=\left\{\begin{array}{ll}
                      1 & \mbox{if \ \ } i\geq r\\
                       0 &\mbox{if \ \ } i<r.
                \end{array}
        \right. \]
We also set $\delta_i(T)=1.$ Let $$O_{g,n}(V)=\< \ \ u\circ_{g,n}
v,\ \ L(-1)u+L(0)u\ \ |\ \ \mbox{homogeneous}\ \ u\in V^r, \
\mbox{and}\ \ v\in V \ \ \>.$$ Here, $u\circ_{g,n}
v=\Res_{z}Y(u,z)v\frac{(1+z)^{\wt u-1+\delta_i(r)+l+r/T}}{z^{2l
+\delta_{i}(r)+\delta_{i}(T-r) }}.$

Define $A_{g,n}(V)=V/O_{g,n}(V).$ Then
$A_{g,n}(V)$ is the untwisted associative algebra $A_n(V)$ as defined
in
[DLM3] if $g$ is the identity element and is $A_g(V)$ in [DLM2] if $n=0.$
We also define a second product $*_{g,n}$ on $V$ for $u$ and $v$ as
above:
\[u*_{g,n}v=\left\{\begin{array}{ll}
        \sum_{m=0}^{l}(-1)^m{m+l\choose
l}\Res_zY(u,z)\frac{(1+z)^{\wt\,u+l}}{z^{l+m+1}}v &\mbox{if}\ \ r=0\\
        0&\mbox{if}\ \ r>0.
\end{array}
\right.\]
Extend this linearly to obtain a bilinear product  on $V.$

Let $M=\oplus_{n\in{1\over T}\Z_{+}}M(n)$ be an admissible $g$-twisted
$V$-module such that $M(0)\neq 0$.
Following
[Z] we define weight zero operator $o_M(v)=v_{\wt v-1}$ on $M$ for
homogeneous $v$ and extend $o_M(v)$ to all $v$ by linearity. It is
clear from the definition that $o_M(v)M(n)\subset M(n)$ for
all $n.$
\begin{th}\label{l2.3} Let $V$ be a vertex operator algebra and
$M$ an admissible $g$-twisted $V$-module. Then

1) The product $*_{g,n}$ induces the structure of an associative
algebra  on $A_{g,n}(V)$ with identity ${\bf 1}+O_{g,n}(V).$
Moreover, $\omega+O_{g,n}(V)$ is a central element of $A_{g,n}(V).$

2) For $0\leq m\leq n,$ the map $\psi_n: v+O_{g,n}(V)\mapsto o_M(v)$ from
$A_{g,n}(V)$ to $\End M(m)$ makes $M(m)$ an $A_{g,n}(V)$-module.

3) $M$ is irreducible $g$-twisted $V$-module if and only if $M(n)$ is
a simple $A_{g,n}(V)$-module
for all $n\in {1\over T}\Z_+.$
\end{th}
Parts 1), and 2) are proved in [DLM4], and 3) is proved in [DM2].

\begin{th}\label{imp}{[DLM4]} There is a bijection map between the isomorphism classes of
the irreducible admissible $g$-twisted
$V$-modules and the isomorphism classes of the simple $A_{g,n}(V)$-modules which can not factor through
$A_{g,n-1/T}(V)$.
\end{th}

\begin{th}\label{thagn}{[DLM4]}
Suppose $V$ is a $g$-rational. Then

1) $A_{g,n}(V)$ is a finite-dimensional semisimple associative for
all $n\in{1\over T}\Z_+$.

2) There is a bijection map between the category of
finite-dimensional $A_{g,n}(V)$-modules whose irreducible
components cannot factor though $A_{g,n-1/T}(V)$ and the category
of ordinary $g$-twisted $V$-modules.

\end{th}

\section{ An application of skew group algebras in orbifold theory}

Let $V$ be a simple vertex operator algebra. Let $G$ be a finite
automorphism group of $V$. In this section we will use the theory
of skew group algebras to show the complete reducibility of any
irreducible $g$-twisted $V$-module as a $V^G$-module when $g$ is
in the center of $G$ and $V$ is a $g$-rational vertex operator
algebra.

We begin with the following setting. For $g\in G$, we let $(M,
Y_M)$ be an irreducible $g$-twisted $V$-module. For $h\in G$, we
set
$$(M,Y_M)\circ h=(M\circ h, Y_{M\circ h}). $$
Here, $M\circ h=M$ as vector spaces and
\begin{equation}\label{sv2}
Y_{M\circ h}(v,z)=Y_M(hv,z)\ \ \mbox{for}\ \ v\in V.
\end{equation}
In fact, the space $M\circ h$ is an irreducible $h^{-1}gh$-twisted $V$-module.
Thus, $M\circ h$ is an irreducible $g$-twisted $V$-module if and only if $h$ is in the centralizer of $g$.

We set $G_M=\{ \ \ h\in G\ \ |\ \ (M,Y_M)\cong (M\circ h,Y_{M
\circ h})\ \ \mbox{as}\ \ g\mbox{-twisted}\ \ V\mbox{-modules}\ \
\}.$ It was proved in [DM3] that $g\in G_M$. To be more precise,
we write $M=\bigoplus_{n=0}^{\infty}M_{\l+{n\over T}}$ where $T$
is the order of $g$. Define $\p(g):M\rightarrow M$ by
$\p(g)|_{M_{\l+{n\over T}}}=e^{-2\pi i {n\over T}}.$ Then $\p(g)$
is an automorphism of $(M,Y_M)$.

For $h\in G_M$, there is a linear isomorphism $\p(h):M\rightarrow M$ satisfying
\begin{equation}\label{es}
\p(h)Y_M(v,z)\p(h)^{-1}=Y_{M\circ h}(v,z)=Y_M(hv,z) \end{equation}
for $v\in V$. The simplicity of $M$ together with Schur's lemma
shows that $h\rightarrow \p(h)$ is a projective representation of
$G_M$ on $M$. Note that we can choose $\p (1)=1_M$. Let $\a_M\in
Z^2(G_M,\C^*)$ be the corresponding 2-cocycle whose values are
roots of unity. Then $M$ is a module for $\C^{\a_M}[G_M]$.
\begin{rem}
Since each $\p(h)$ commutes with $L(n)$ for all $n\in \Z$, $\p(h)$ preserves the homogeneous subspaces $M(m)$, and $M(m)$ is a $\C^{\a_M}[G_M]$-module.
\end{rem}
\begin{rem}\ \

1) If $M=V$, we can take $\p(h)=h$ and $\C^{\a_M}[G_M]=\C[G].$

2) Let $G=\<h\>$ be a cyclic group of prime order $p$. For $1\leq i\leq p-1$, we let $M$ be an irreducible $h^i$-twisted $V$-module. Then $G_M=\<h\>$ and the 2-cocycle $\a_M$ can be taken to be trivial.
\end{rem}

{\em Assume that $g$ is an element in the center of $G$} and {\em
$V$ is $g$-rational}. Let $M$ be an irreducible $g$-twisted
$V$-module. Then $M(n)$ is a simple $A_{g,n}(V)$ for all $n\in
{1\over T}\Z_+$. Moreover, by equations (\ref{sk2}), (\ref{sv2}),
we have $$^hM(n)=M\circ {h^{-1}}(n)\ \ \mbox{for}\ \ h\in
G,n\in{1\over T}\Z_+.$$ We let $$G_{M(n)}=\{ \ \ h\in G \ \ |\ \
^hM(n)\cong M(n)\ \ \mbox{as} \ \ A_{g,n}(V) \mbox{-modules}\ \
\}$$ be the inertia subgroup of $M(n)$.
\begin{lem} For $n\in
{1\over T}\Z_+$, if $M(n)\neq 0$, we have $G_{M(n)}=G_M.$
\end{lem}
\pf Recall that $M$ and $M\circ h$ are irreducible $g$-twisted
$V$-modules. By Theorem \ref{imp}, we conclude that $M\cong M\circ
h$ as $g$-twisted $V$-modules if and only if $M(n)\cong M\circ
h(n)$ as $A_{g,n}(V)$-modules. Thus $G_{M(n)}=G_M$ if $M(n)\neq 0$.\qed

Let $\L_{G_M,\a_M}$ be the set of all irreducible characters $\l$
of $\C^{\a_M}[G_M]$. For each $\l\in\L_{G_M,\a_M}$, we denote the
corresponding simple module for $\l$ by $W_{\l}$ and we let
$M^{\l}$ be the sum of simple $\C^{\a_M}[G_M]$-submodules of $M$
isomorphic to $W_{\l}$. Since $M^{\l}$ is nonzero for any $\l\in
\L_{G_M,\a_M}$ (cf. [DY]), we see that
$$M=\bigoplus_{\l\in \L_{G_M,\a_M}}M^{\l}.$$
As in section 3, we let $M_{\l}=\{f(w) |f\in \Hom_{\C^{\a_M}[G_M]}(W_{\l},M)\}$ for a fixed nonzero $w\in W_{\l}$.

\begin{th} For any $\g\in \L_{G_M,\a_M}$, $M_{\g}$ is a $V^{G}$-module.
\end{th}

\pf By equation (\ref{es}), we have $\p(h)Y_M(v,z)=Y_M(v,z)\p(h)$
for all $v\in V^G$. Thus the action of $\C^{\a_M}[G_M]$ and $V^G$
are commute on $M$. Hence, $M_{\g}$ is a $V^G$-module. \qed

For any $\l\in \L_{G_M,\a_M}$, we identify $M^{\l}$ with $M_{\l}\o W_{\l}$ as $V^G\o \C^{\a_M}[G_M]$-modules.
Thus $$M=\bigoplus_{\l\in\L_{G_M,\a_M}}W_{\l}\o M_{\l}.$$
\begin{th} For any $\g\in \L_{G_M,\a_M}$, $M_{\g}$ is an irreducible
$V^{G}$-module.

\end{th}
\pf
For $n\in{1\over T}\Z_+$, we set $M^{\g}(n)=M^{\g}\bigcap M(n)$ and
$M_{\g}(n)=M_{\g}\bigcap M(n)$.
In order to show that $M_{\g}$ is an irreducible $V^{G}$-module, it
is enough to show that for any $n\in {1\over T}\Z_+$, if $M^{\g}(n)$
is nonzero then $M_{\g}(n)$ is a simple $A_{g,n}(V^G)$-module (cf. Theorem \ref{l2.3} 3) ).
We now take $n\in {1\over T}\Z_+$ such that $M^{\g}(n)\neq 0$.
Since $V$ is $g$-rational, $A_{g,n}(V)$ is a finite-dimensional
semisimple associative algebra (cf. Theorem \ref{thagn}).
Since $M(n)$ is a simple $A_{g,n}(V)$-module, Theorem \ref{thms} implies that $M_{\g}(n)$ is a simple
$A_{g,n}(V)^G$-module.
Since the identity on $V^G$ induces an algebra epimorphism
from $A_{g,n}(V^G)$ to $A_{g,n}(V)^G$, we can conclude $M_{\g}(n)$ is a simple $A_{g,n}(V^G)$-module.
Consequently, $M_{\g}$ is an irreducible $V^G$-module. \qed

\begin{cor}
$M$ is a completely reducible $V^G$-module.
\end{cor}
\begin{cor}
If $V$ is a rational vertex operator algebra then every irreducible $V$-module is a completely reducible $V^G$-module.
\end{cor}
\begin{cor}\ \

1) Assume $g$ is in the center of $G$. Let $M$ be an irreducible
$g$-twisted $V$-module. If $G_M$ is a cyclic group, says $\<h\>$,
then, in fact, the eigenspaces of $M$ with respect to the action
of $h$ are irreducible $V^G$-modules.

2) Suppose $G$ is a cyclic group of prime order $p$. Let $M$ be an irreducible $g^i$-twisted $V$-module where $1\leq i\leq p-1$. Hence the eigenspaces of $M$ are irreducible $V^G$-modules.
\end{cor}

\end{document}